\documentclass[12pt]{amsart}
\usepackage[centertags]{amsmath}
\usepackage{amsfonts,amssymb}
\usepackage{graphicx}
\usepackage{enumerate}
\usepackage[latin1]{inputenc}

  \newtheorem{theorem}{Theorem}
  \newtheorem{corollary}{Corollary}
  \newtheorem{proposition}{Proposition}
  \newtheorem{lemma}{Lemma}%
  \theoremstyle{remark}
  \newtheorem{remark}{Remark}

  \newcommand{\Sc}{\mathcal S}
  \newcommand{\Tc}{\mathcal T}
  \newcommand{\Ac}{\mathcal A}
  \newcommand{\Bc}{\mathcal B}
  
  \newcommand{\Dc}{\mathcal D}
  \newcommand{\Fc}{\mathcal F}
  \newcommand{\Pc}{\mathcal P}
 \newcommand{\Uc}{\mathcal U}
\newcommand{\Cc}{\mathcal C}
\begin{document}

\title{The denominators of the Bernoulli numbers}
\date{\today}

\author{Carl Pomerance}
\address{Mathematics Department, Dartmouth College, Hanover, NH 03784}
\email{carl.pomerance@dartmouth.edu}

\author{Samuel S. Wagstaff, Jr.}
\address{Center for Education and Research in Information Assurance and Security
  and Department of Computer Sciences, Purdue University \\
  West Lafayette, IN 47907-1398 USA}
\thanks{S.S.W.'s work was supported by the CERIAS Center at Purdue University}

\email{ssw@cerias.purdue.edu}

\begin{abstract}
We study the asymptotic density of the set of subscripts of the
Bernoulli numbers having a given denominator.  We also study the distribution
of distinct Bernoulli denominators and some related problems.
\end{abstract} 
\maketitle

\section{Introduction}

Long before the proof of Andrew Wiles, it was thought that the
path to Fermat's Last Theorem (FLT) led through the Bernoulli numbers.
  Defined by the series
\[
\frac t{e^t-1}=\sum_{n=0}^\infty B_n\frac{t^n}{n!},
\]
the Bernoulli numbers $B_n$ are rationals, in lowest terms $N_n/D_n$,
and both the sequence of numerators $N_n$ and denominators $D_n$ have a 
connection to FLT.  It is known since Ernst Kummer, that FLT holds for any
odd prime that does not divide the class number of the cyclotomic
field ${\mathbb Q}[e^{2\pi i/p}]$ (such primes are called ``regular"),
and that a condition for this to occur is that $p\nmid N_n$ for all even
$n\le p-3$.  The so-called first case of FLT for a prime exponent $p$
is when $p$ divides none of the powers.  Sophie Germain proved this
case when $2p+1$ is also a prime; we now call such primes $p$
Germain primes.  If $p$ is a Germain prime, then $2p+1\mid D_{2p}$,
giving the connection to Bernoulli denominators.

This paper considers the density of $n$ with a fixed value of $D_n$, 
the distribution of distinct values of $D_n$, and some other related
problems.  Note that $B_0=1$, $B_1=-1/2$,
and $B_{2k+1}=0$ for all integers $k>0$.  So, we only consider the remaining
cases $D_{2k}$, $k>0$.  We have a
precise formula for $D_{2k}$ in these cases given by the theorem
of Karl von Staudt and Thomas Clausen: {\it $D_{2k}$ is
the product of the primes $p$ for which $p-1\mid 2k$.}  

 From the von Staudt--Clausen theorem we immediately see that $D_{2k}$ is
 squarefree and a multiple of 6.  We set some notation.
For a positive integer $m$ let $T_m = \{p\ {\mathrm{prime}} : p-1 \mid m\}$.
Thus, $T_m = \{2\}$ when $m$ is odd.  For even $m$, we have, for example,
$T_2 = \{2,3\}$,
$T_4 = \{2,3,5\}$, $T_6 = \{2,3,7\}$, $T_8 = T_4$, $T_{10} = \{2,3,11\}$, etc.
As we have seen, $D_{2k}$ is the product of the primes in $T_{2k}$.

For $n>0$ even, let
\[
\Sc_{n} = \{m > 0\hbox{ even} : T_{m} = T_{n}\}=\{m>0\hbox{ even}:D_{m}=D_{n}\}.
\]
Then $\Sc_2 = \{2,14,26,34,38,\ldots\}$, $\Sc_4 = \{4,8,68\ldots\}$, etc.
So $\Sc_{n}$ is the set of all even $m$ for which $B_{m}$ has the same
denominator as $B_{n}$, namely, the product of the primes in $T_{n}$.
In 1980, Erd\H os and the second author proved \cite{EW}
that $\Sc_{n}$ has a positive asymptotic density for every even $n$.
Sunseri \cite{Sun} proved that the density
of $\Sc_2$ is at least as large as the density of $\Sc_n$
for every even $n > 0$.
We will give a simple proof below of a slightly stronger version of Sunseri's result.

Let 
\[
\Dc:=\{D_n:n>0~\hbox{even}\}=\{6,30,42,66,\dots\}.
\]
Further, for $d\in\Dc$, let $F_d=\min\{n:D_n=d\}$.  For example,
$F_6=2$, $F_{30}=4$, $F_{42}=6$, and $F_{66}=10$.  Let
\[
\Fc:=\{F_d:d\in\Dc\}=\{2,4,6,10,\dots\}.
\]
Let $\lambda$ denote the Carmichael $\lambda$-function.  In particular,
for $n$ squarefree, $\lambda(n)=\hbox{lcm}\{p-1:p\mid n\}$, where $p$ denotes
a prime variable.  We characterize $\Fc$ as the set of values of $\lambda(n)$
for $n>2$ squarefree, and use this characterization plus some results on the
distribution of the image of $\lambda$ to get good approximations to the
counting functions of $\Dc$ and $\Fc$.

In this paper $p$ always denotes a prime.
For $p>2$, we have $D_{p-1}=dp$ for some even integer $d$.
We show that but for a set of primes of relative density 0 in the set of primes,
this number $d$ itself is in $\Dc$.  Further, we show that for each fixed $d\in\Dc$,
the relative density of the primes $p$ with $D_{p-1}=dp$ is positive,
and the sum of these densities is 1.

\section{Characterization of $\Fc$}

Here we give a connection between the set $\Fc$ and the image of the
Carmichael $\lambda$-function.

\begin{proposition}
\label{prop:lambda}
For each $d\in\Dc$ we have $F_d=\lambda(d)$.  Further,
\[
\Fc=\{\lambda(n):n>2\hbox{ squarefree}\}.
\]
\end{proposition}
\begin{proof}
Let $d\in\Dc$ and suppose $n$ has $D_n=d$.  If $p$ is a prime factor of $d$,
then we have $p-1\mid n$, by von Staudt--Clausen.  
Thus, $\lambda(d)\mid n$.  Clearly if $a\mid b$,
with $a,b$ even, then $D_a\mid D_b$.  Thus, $D_{\lambda(d)}\mid D_n=d$.
Also, $p\mid d$ implies that $p-1\mid\lambda(d)$, which implies that
$p\mid D_{\lambda(d)}$.  Since $d$ is squarefree, we thus have 
$d\mid D_{\lambda(d)}$.  Hence, $d=D_{\lambda(d)}$
and $\lambda(d)=F_d$.  This proves the first assertion and half of the second
assertion.  Say $m>2$ is squarefree and $n=\lambda(m)$.  Let $d=D_n$.
By the von Staudt--Clausen theorem, $m\mid d$, so that 
$n=\lambda(m)\mid\lambda(d)$.  As we saw above, whenever we have
$d=D_n$, we have $\lambda(d)\mid n$.  Thus, $n=\lambda(d)=\lambda(m)$
and the proposition follows.
\end{proof}


\section{$\Sc_2$ has the greatest density}

We introduce more notation.  Let $\Ac$ denote a set of positive integers and write
$\Bc(\Ac)$ for the set of all positive integer
multiples of the elements of $\Ac$.
Let $A(x)$ denote the counting function of the set $\Ac$.
Write ${\rm d}(\Ac)=\lim_{x\rightarrow\infty} A(x)/x$ for
the asymptotic density of $\Ac$, if it exists.   For a real number $t$,
write $\Ac^{(t)}$ for $\{a\in\Ac:a>t\}$.

We record the following result from \cite[Lemmas 1, 2]{PS}.
\begin{proposition}
\label{prop:PS}
Suppose that $\Ac$ is a set of positive integers not containing $1$ with
the property that 
\begin{equation}
\label{eq:larget}
\lim_{t\to\infty}\limsup_{x\to\infty}\Bc(\Ac^{(t)})(x)/x=0.
\end{equation}
Then ${\rm d}(\Bc(\Ac))$ exists and is $<1$.
\end{proposition}
Note that the condition \eqref{eq:larget} essentially
asserts that most numbers are not divisible by a large member of $\Ac$.

 If $n$ is even, let 
 \[
 \Uc_n=\{mn:D_{mn}=D_n\}.
 \]
 Note that if $n\in\Fc$, then $\Uc_n=\Sc_n$.  
For $r$ prime, let
  \[
 \Uc_{n,r}=\{mn\in\Uc_n:r\,\nmid\,m\}.
 \]
 \begin{lemma}
 \label{lem:dens1}
 With the above notation, the sets $\Uc_n$ and $\Uc_{n,r}$ have positive
 asymptotic densities.
 \end{lemma}
\begin{remark}
When $n\in\Fc$, the result that $\Uc_n$ has positive asymptotic density
follows from  \cite[Theorem 3]{EW}.
 \end{remark}
\begin{proof}
Let 
\[
\Ac_n=\{(p-1)/\gcd(p-1,n):p-1\nmid n\}.
\]
It follows from \cite[Theorem 2]{EW} that condition \eqref{eq:larget}
holds for $\Ac=\Ac_n$.   Indeed, $\Ac_n$ is the disjoint union over
the divisors $d$ of $n$ of the sets
$\Ac_{n,d}:=\{(p-1)/d:d=\gcd(p-1,n),\,p-1\nmid n\}$, and since \eqref{eq:larget}
holds for each separate $\Ac_{n,d}$ by \cite[Theorem 2]{EW}, and
since $n$ is fixed, it holds too for $\Ac_n$.
Thus, from Proposition \ref{prop:PS} we have that $\Bc(\Ac_n)$ has an
asymptotic density $<1$.  Let $\Cc_n={\mathbb N}\setminus\Bc(\Ac_n)$.
Then $\Cc_n$ has a positive asymptotic density.  Our first assertion now
follows upon noting that $\Uc_n=n\Cc_n=\{nm:m\in\Cc_n\}$.  Indeed, if
$p-1\mid n$ then $p-1\mid nm$ for all $m$, and if $p-1\mid nm$ for some
$p$ with $p-1\,\nmid\, n$, then $(p-1)/\gcd(p-1,n)\mid m$, so that $m\notin\Cc_n$.

The second assertion follows from the same argument applied to
$\Ac=\Ac_n\cup\{r\}$.
\end{proof}
 
 \begin{lemma}
 \label{lem:trythis}
 If $r$ is prime and $n$ is even,  then
 \[
 {\rm d}(\Uc_{nr,r})\le\frac1r{\rm d}(\Uc_{n,r})
~ \hbox{ and }~{\rm d}(\Uc_n)\le\frac r{r-1}{\rm d}(\Uc_{n,r}).
 \]
 Equality in the second assertion occurs only if each $D_{nr^i}=D_n$
 for all $i\ge0$.
 \end{lemma}
 \begin{remark}
 It is not clear if there are any pairs $n,r$ with $D_{nr^i}=D_n$ for all $i\ge0$.
 However, if there are only finitely many Fermat primes, with $p=2^{2^k}+1$
 being the largest one, then $D_{(p-1)2^i}=D_{p-1}$ for all $i\ge0$.
 \end{remark}
 \begin{proof}
 Consider the map $m\mapsto m/r$ on $\Uc_{nr,r}$.  Let $m\in\Uc_{nr,r}$.
  We'd like to show that $m/r\in\Uc_{n,r}$.  Suppose $r^J\,\|\,n$.
 Clearly $r^J\,\|\,m/r$ and $n\mid m/r$,
  so it remains to show that $D_{m/r}=D_{n}$.  
 If $p-1\mid n$ then $p-1\mid m/r$, since $n\mid m/r$.  Further, if
 $p-1\mid m/r$, then $p-1\mid m$, so that $p-1\mid nr$, since $D_{m}=D_{nr}$.  
 But $r^{J+1}\nmid m/r$,
 so that $r^{J+1}\nmid p-1$.  That is, $p-1=ur^i$, where $u\mid n/r^J$ and $i\le J$.
 Hence $p-1\mid n$ and $m/r\in\Uc_{n,r}$.  The
 first assertion of the lemma follows.
 
 For the second, note that $\Uc_n$ is contained in the disjoint union
 of those sets $\Uc_{nr^i,r}$ for $i=0,1,\dots$ with $D_{nr^i}=D_n$.  
 By the first part of
 the lemma applied $i$ times, ${\rm d}(\Uc_{nr^i,r})\le r^{-i}{\rm d}(\Uc_n)$.
 It remains to note that $\sum_{i\ge0}r^{-i}=r/(r-1)$.  Note that in the possible
 case that each $D_{nr^i}=D_n$ we are expressing $\Uc_n$ as an infinite
 disjoint union of sets with positive asymptotic density.  Asymptotic density
 is not necessarily countably additive, but in this case there is no issue since
$\cup_{i\ge k}\,\Uc_{nr^i,r}$ is contained in the multiples of $r^k$ and so 
has upper density which tends to 0 as $k\to\infty$.
 \end{proof}
 
 \begin{theorem}
\label{thm:ab}
If $a,b\in\Fc$  with $a\mid b$, then 
\begin{equation}
\label{eq:phi}
{\rm d}(\Sc_b)\le\frac1{\varphi(b/a)}{\rm d}(\Sc_a).
\end{equation}
In addition, ${\rm d}(\Sc_4)\le\frac34{\rm d}(\Sc_2)$.
\end{theorem}
\begin{remark}
Since every member of $\Fc$ is even, Theorem \ref{thm:ab} in the case $a=2$
shows that ${\rm d}(\Sc_2)$ is at least one-third larger than every
other ${\rm d}(\Sc_n)$.
\end{remark}
\begin{proof}
Let $r$ be a prime factor of $b/a$ and let $K=v_r(b/a)$.
By the second part of Lemma \ref{lem:trythis},
\[
{\rm d}(\Sc_b)={\rm d}(\Uc_b)\le\frac r{r-1}{\rm d}(\Uc_{b,r}).
\]
  Repeatedly
applying the first part of Lemma \ref{lem:trythis},
\[
{\rm d}(\Uc_{b,r})\le\frac1r{\rm d}(\Uc_{b/r,r})\le\dots
\le \frac1{r^{K}}{\rm d}(\Uc_{b/r^{K},r}),
\]
so that
\[
{\rm d}(\Uc_b)\le\frac1{\varphi(r^{K})}{\rm d}(\Uc_{b/r^{K},r}).
\]
Thus, \eqref{eq:phi} follows by induction on the number of distinct prime divisors of $b/a$.

Now suppose that $a=2,b=4$, and follow the above proof.  Note that
$D_{16}$ is divisible by 17, but $D_4$ is not.  Thus, $\Sc_4=\Uc_{4,2}\cup\,\Uc_{8,2}$, and 
\[
{\rm d}(\Sc_4)={\rm d}(\Uc_{4,2}\cup\,\Uc_{8,2})\le\frac32{\rm d}(\Uc_{4,2}).
\]
Also
\[
{\rm d}(\Uc_{4,2})\le\frac12{\rm d}(\Uc_{2,2})=\frac12{\rm d}(\Uc_2)=\frac12{\rm d}(\Sc_2).
\]
Hence ${\rm d}(\Sc_4)\le\frac34{\rm d}(\Sc_2)$.
\end{proof}

\begin{remark}
The proof shows that ${\rm d}(\Sc_6)\le\frac12{\rm d}(\Sc_2)$ and
${\rm d}(\Sc_{10})\le\frac14{\rm d}(\Sc_2)$, for example,
but provides no way to compare ${\rm d}(\Sc_4)$ with 
${\rm d}(\Sc_6)$ or ${\rm d}(\Sc_6)$ with ${\rm d}(\Sc_{10})$.
\end{remark}

\begin{corollary}
\label{cor:denom6}
Measured by the asymptotic density of their sets of subscripts,
more Bernoulli numbers have denominator $6$ than any other integer.
\end{corollary}

\begin{remark}
\label{rmk:densum}
Every even number $n$ is in some $\Sc_f$ for $f\in\Fc$, namely for
$f=\lambda(D_n)$.  Moreover, the sets $\Sc_f$ for $f\in\Fc$ are
pairwise disjoint.  It follows that $\sum_{f\in\Fc}{\rm d}(\Sc_f)\le\frac12$.
In fact, we have
\begin{equation}
\label{eq:sum}
\sum_{f\in\Fc}{\rm d}(\Sc_f)=\frac12,
\end{equation}
 as asserted in \cite[Corollary, p.\ 111]{EW}.
The proof follows immediately from \cite[Theorem 2]{EW}, which asserts\footnote
{This result has been sharpened in the recent papers \cite{F2}, \cite{MPP}.}
that for each $\epsilon>0$ there is some number $B$ such that the upper density
of those integers divisible by some $p-1>B$ is $<\epsilon$.  (That is, \eqref{eq:larget} 
holds.)  So, to prove
\eqref{eq:sum}, note that there are only finitely
many $d\in\Dc$ not divisible by any prime $p>B+1$, since members of $\Dc$
are squarefree, so there are only finitely many $f\in\Fc$ not divisible by any
$p-1>B$.  The numbers $n$ in all other sets $\Sc_f$ are divisible by some
$p-1>B$ so they comprise a set of upper density $<\epsilon$.  Thus, the
sum in \eqref{eq:sum} is $>\frac12-\epsilon$.
\end{remark}

\section{Statistics for the $\Sc_{n}$}

If one defines $\Sc_1$ as $\{n:T_n=T_1=\{2\}\}$, then we merely have
the set of odd numbers, having 
density $\frac12$.
What about the densities of the various $\Sc_n$ with $n$ even?  We have
proved that $\Sc_2$ has density at least $\frac13$ more than the next largest
density.
Numerical calculation suggests that $\Sc_2$ has density about $0.07$,
followed by $\Sc_4$ (0.03) and $\Sc_6$ (0.01) in that order.
We know from the proof of Theorem \ref{thm:ab} that ${\rm d}(\Sc_4)$ is
the largest of the ${\rm d}(\Sc_n)$ for $4\mid n$ and that ${\rm d}(\Sc_6)$
is the largest of the ${\rm d}(\Sc_n)$ for $6\mid n$.  As mentioned, the proof does not
allow us to compare ${\rm d}(\Sc_4)$ and ${\rm d}(\Sc_6)$.
It follows from the methods in \cite{EW}, \cite{Sun} that the densities
are effectively computable in principle, and so in principle it is possible
to resolve this issue, though such calculations appear currently to be infeasible.

Tables \ref{Ta:count579} and \ref{Ta:count579a} give the number of elements
in $\Sc_{n}$ for $n=2k\le112$ up to $10^m$ for $m=5$, 7 and 9.
In these tables, $2k$ and $2k'$ are the smallest two elements of $\Sc_{2k}$,
so that $2k\in\Fc$.
The partial densities may be computed easily from the counts.
Note that, in the range from $10^5$ to $10^9$, the partial
densities tend to decrease when $T_{2k}$ contains few primes
and increase when $T_{2k}$ contains many primes.

The tables suggest that
\[ 
{\rm d}(\Sc_2)> {\rm d}(\Sc_4)> {\rm d}(\Sc_6)>
{\rm d}(\Sc_{10})> {\rm d}(\Sc_{16})
\]
are the five largest densities of the $\Sc_{2k}$.
In other words, the five most popular denominators of Bernoulli
numbers seem to be 6, 30, 42, 66, 510, in that order.  

\begin{table}[ht]
\caption{Number of elements of $\Sc_{2k}$ below various limits.}
\label{Ta:count579}
\begin{tabular}{rrcrrr}
First & Second & $T_{2k}$ & Count & Count & Count \\
$2k$ & $2k'$ & & $\le10^5$ & $\le10^7$ & $\le10^9$ \\ 
2 & 14 & $\{2,3\}$ & 7992 & 758582 & 73129588 \\
4 & 8 & $\{2,3,5\}$ & 3423 & 320500 & 30579077 \\
6 & 114 & $\{2,3,7\}$ & 1371 & 125712 & 11923816 \\
10 & 50 & $\{2,3,11\}$ & 1080 & 99675 & 9457553 \\
12 & 24 & $\{2,3,5,7,13\}$ & 495 & 49498 & 4751091 \\
16 & 32 & $\{2,3,5,17\}$ & 713 & 67742 & 6379485 \\
18 & 54 & $\{2,3,7,19\}$ & 397 & 38502 & 3671790 \\
20 & 340 & $\{2,3,5,11\}$ & 289 & 27745 & 2609924 \\
22 & 154 & $\{2,3,23\}$ & 566 & 52508 & 4959735 \\
28 & 56 & $\{2,3,5,29\}$ & 309 & 29692 & 2793858 \\
30 & 1770 & $\{2,3,7,11,31\}$ & 138 & 13615 & 1309849 \\
36 & 3924 & $\{2,3,5,7,13,19,37\}$ & 72 & 7846 & 799642 \\
40 & 6680 & $\{2,3,5,11,41\}$ & 92 & 10044 & 950144 \\
42 & 294 & $\{2,3,7,43\}$ & 124 & 12645 & 1199553 \\
44 & 484 & $\{2,3,5,23\}$ & 160 & 15325 & 1433972 \\
46 & 322 & $\{2,3,47\}$ & 261 & 24295 & 2290634 \\
48 & 10128 & $\{2,3,5,7,13,17\}$ & 26 & 4572 & 497209 \\
52 & 104 & $\{2,3,5,53\}$ & 164 & 16638 & 1558130 \\
58 & 406 & $\{2,3,59\}$ & 235 & 20607 & 1935087 \\
\end{tabular}
\end{table}

\begin{table}[ht]
\caption{Number of elements of $\Sc_{2k}$ below various limits.}
\label{Ta:count579a}
\begin{tabular}{rrcrrr}
First & Second & $T_{2k}$ & Count & Count & Count \\
$2k$ & $2k'$ & & $\le10^5$ & $\le10^7$ & $\le10^9$ \\ 
60 & 13620 & $\{2,3,5,7,11,13,31,61\}$ & 21 & 2917 & 340111 \\
66 & 3894 & $\{2,3,7,23,67\}$ & 77 & 7202 & 680301 \\
70 & 350 & $\{2,3,11,71\}$ & 83 & 8815 & 818849 \\
72 & 12024 & $\{2,3,5,7,13,19,37,73\}$ & 12 & 2137 & 259257 \\
78 & 1014 & $\{2,3,7,79\}$ & 71 & 6771 & 636574 \\
80 & 160 & $\{2,3,5,11,17,41\}$ & 39 & 5960 & 610485 \\
82 & 574 & $\{2,3,83\}$ & 150 & 13715 & 1293383 \\
84 & 168 & $\{2,3,5,7,13,29,43\}$ & 16 & 2924 & 339634 \\
88 & 968 & $\{2,3,5,23,89\}$ & 53 & 5593 & 528007 \\
90 & 14670 & $\{2,3,7,11,19,31\}$ & 17 & 2629 & 284131 \\
92 & 184 & $\{2,3,5,47\}$ & 116 & 10822 & 1017455 \\
96 & 20256 & $\{2,3,5,7,13,17,97\}$ & 7 & 1645 & 196489 \\
100 & 1700 & $\{2,3,5,11,101\}$ & 34 & 4115 & 393270 \\
102 & 1734 & $\{2,3,7,103\}$ & 50 & 5041 & 473949 \\
106 & 1378 & $\{2,3,107\}$ & 120 & 10794 & 1007709 \\
108 & 11772 & $\{2,3,5,7,13,19,37,109\}$ & 14 & 1593 & 190046 \\
110 & 550 & $\{2,3,11,23\}$ & 72 & 6481 & 609261 \\
112 & 224 & $\{2,3,5,17,29,113\}$ & 41 & 4135 & 422188 \\
\end{tabular}
\end{table}

To form $\Sc_2$, begin with the set $\Tc$ of even positive integers.
This is the set of positive integers divisible by both $2-1$ and $3-1$.
Now delete from $\Tc$ all multiples of $q-1$ for all primes $q>4$.
Some primes $q>4$ may be skipped because the multiples of $q-1$ 
were deleted when multiples of $r-1$ were removed for some prime $r<q$.
For example, multiples of $13-1$ were deleted when multiples of $7-1$
were removed.  However, we must always delete the multiples of $q-1$ 
whenever $p=(q-1)/2$ is prime.
Recall that the prime $p$ is called a Germain prime if $2p+1$ is also prime.
Thus, it is necessary to remove all multiples of the Germain primes from
$\Tc$ to form $\Sc_2$.
But it is not sufficient to delete multiples of Germain
primes because, for example, $q=239$ is prime
but $p=(q-1)/2=119$ is not prime and not divisible by a Germain prime.

Table \ref{Ta:countmod} lists the number of elements of $\Sc_{2k}$
below $10^6$ in each residue class modulo 8 and the first few odd primes.
(The residue classes $>7$ modulo 11 and 13 are omitted to make the table
fit on a page. The missing values are similar to the last ones shown
in that row of the table.)
The value ``total'' is the number of elements of $\Sc_{2k}$ less than $10^6$.

\begin{footnotesize}
\begin{table}[ht]
\caption{Number of elements $\le10^6$ of $\Sc_{2k}$ in residue classes.}
\label{Ta:countmod}
\begin{tabular}{|rrrrrrrrr|} \hline
Modulus & 0 & 1 & 2 & 3 & 4 & 5 & 6 & 7 \\ \hline
& \multicolumn{8}{c|}{ $2k = 2$, total = 77696}  \\
8 & 0 & 0 & 38849 & 0 & 0 & 0 & 38847 & 0 \\
3 & 0 & 31612 & 46084 & & & & & \\
5 & 0 & 18636 & 18565 & 19097 & 21398 & & & \\
7 & 9179 & 11168 & 11175 & 11080 & 11309 & 11125 & 12660 & \\
11 & 0 & 7661 & 7671 & 7682 & 7730 & 7726 & 7649 & 7627 \\
13 & 5116 & 5959 & 5980 & 5975 & 5970 & 5972 & 6079 & 6035 \\ \hline
& \multicolumn{8}{c|}{ $2k = 4$, total = 33001 } \\
8 & 9490 & 0 & 0 & 0 & 23511 & 0 & 0 & 0 \\
3 & 0 & 15877 & 17124 & & & & & \\
5 & 0 & 7868 & 7244 & 9186 & 8703 & & & \\
7 & 0 & 5186 & 5176 & 5328 & 5274 & 6089 & 5948 & \\
11 & 0 & 3160 & 3198 & 3206 & 3179 & 3191 & 3200 & 3338 \\
13 & 0 & 2693 & 2633 & 2679 & 2637 & 2695 & 2682 & 2669 \\ \hline
& \multicolumn{8}{c|}{ $2k = 6$, total = 12996 } \\
8 & 0 & 0 & 6508 & 0 & 0 & 0 & 6488 & 0 \\
3 & 12996 & 0 & 0 & & & & & \\
5 & 0 & 2859 & 3346 & 2867 & 3924 & & & \\
7 & 0 & 2012 & 1990 & 1940 & 2351 & 2004 & 2699 & \\
11 & 0 & 1264 & 1243 & 1252 & 1257 & 1239 & 1227 & 1224 \\
13 & 0 & 1042 & 1044 & 1023 & 1025 & 1018 & 1028 & 1068 \\ \hline
& \multicolumn{8}{c|}{ $2k = 10$, total = 10339 } \\
8 & 0 & 0 & 5175 & 0 & 0 & 0 & 5164 & 0 \\
3 & 0 & 4954 & 5385 & & & & & \\
5 & 10339 & 0 & 0 & 0 & 0 & & & \\
7 & 0 & 1545 & 1948 & 1632 & 1643 & 1571 & 2000 & \\
11 & 0 & 987 & 973 & 1009 & 989 & 991 & 1203 & 999 \\
13 & 0 & 846 & 856 & 850 & 836 & 816 & 841 & 830 \\ \hline
\end{tabular}
\end{table}
\end{footnotesize}

The elements of $\Sc_2$ are all $\equiv2\pmod{4}$ and
seem to be equally distributed between $2\pmod{8}$ and $6\pmod{8}$.
For each Germain prime $p$, the residue class $0\pmod{p}$ is empty
because these classes were removed from $\Tc$ when $\Sc_2$ was formed.
The elements of $\Sc_2$ appear to be equally distributed among
other residue classes modulo Germain primes.

The elements of $\Sc_4$ are all divisible by 4 because $5\in T_4$.
At first it was puzzling why  there are so many more of them that are
$4\pmod{8}$ than $0\pmod{8}$.  However, those members of $\Sc_4$
that are $4\pmod8$ comprise $\Uc_{4,2}$ and those that are $0\pmod8$
comprise $\Uc_{8,2}$, using the notation introduced in the previous section.
In the proof of Theorem~\ref{thm:ab} we saw that ${\rm d}(\Uc_{8,2})\le\frac12
{\rm d}(\Uc_{4,2})$, since dividing a member of $\Uc_{8,2}$ by 2 gives a
member of $\Uc_{4,2}$.  This explains part of the favoring of $4\pmod8$ over
$0\pmod8$, but not all.  In fact, multiplying a member of $\Uc_{4,2}$ by
2 does not always give a member of $\Uc_{8,2}$.  
Another force at work here is that if $m>3$ is odd, then $4m\in\Sc_4$ if and only if 
for each $d\mid m$ with $d>3$, both $2d+1$
and $4d+1$ are composite.  But for $8m\in\Sc_4$, there is the additional
requirement that $8d+1$ is composite.

As we mentioned earlier, every element of $\Sc_{2k}$ is divisible by
its least element $2k$.
There are no multiples of 7 in $\Sc_4$, $\Sc_6$ or $\Sc_{10}$
because 29, 43 and 71, respectively, are prime.
Similarly, there are no multiples of 13 in $\Sc_4$, $\Sc_6$ or $\Sc_{10}$
because 53, 79 and 131 are prime.

\section{The distribution of distinct Bernoulli denominators}

Let $D(x), F(x)$
be the counting functions of $\Dc$ and $\Fc$, respectively,
and let
\[
\beta=1-(1+\log\log2)/\log2=0.08607\dots,
\]
 the Erd\H os--Tenenbaum--Ford constant.
\begin{theorem}
\label{thm:Ddist}
We have, as $x\to\infty$,
\[
D(x)=x/(\log x)^{1+o(1)},\quad F(x)=x/(\log x)^{\beta+o(1)}.
\]
In particular, ${\rm d}(\Dc)={\rm d}(\Fc)=0$.
\end{theorem}
\begin{proof}
Note that Proposition \ref{prop:lambda} implies that the function sending $d\in\Dc$ to 
$\lambda(d)\in\Fc$ is a bijection.  Thus,
\[
D(x)=\#\{\lambda(d):d\in\Dc,\,d\le x\}
\le\#\{\lambda(n):n\le x\}.
\]
Further, the second part of Proposition \ref{prop:lambda} implies that
\begin{align*}
F(x)&=\#\{\lambda(n):n>2\hbox{ squarefree},~\lambda(n)\le x\}\\
&\le\#\{\lambda(n):n\hbox{ such that }\lambda(n)\le x\}.
\end{align*}
 In \cite[Theorem 1.3]{LP} it is shown that 
 \[
 \#\{\lambda(n):n\le x\}=x/(\log x)^{1+o(1)}\hbox{ as }x\to\infty,
 \]
 and \cite[Theorem 1.1]{LP} implies that
\begin{equation}
\label{eq:flp}
\#\{\lambda(n):\lambda(n)\le x\}\le x/(\log x)^{\beta+o(1)}\hbox{ as }x\to\infty.
\end{equation}
So the upper bounds implicit in the theorem follow.
In \cite{FLP} it is shown that equality holds in \eqref{eq:flp}.   In fact, the
method of proof gives the same bound for $\#\{\lambda(n)\le x:n\hbox{ squarefree}\}$,
so by this result, the proof for $\Fc$ is complete.
We show in Theorem \ref{thm:dens} below that $\Dc$ contains the numbers $6p$ 
for a positive proportion of the primes $p$, so the lower bound for $\Dc$
will follow from the prime number theorem.
\end{proof}

\begin{remark}
Looking at small values of $\Fc$ it seems that many are of the form $p-1$
with $p$ prime.  Every $p-1$ is in $\Fc$ for $p>2$ prime, as is easily seen,
but Theorem \ref{thm:Ddist} shows that most members of $\Fc$ are not
in this form. 
\end{remark}

Table \ref{Ta:countF} shows the growth rate of $\Fc$,
the set of first elements $2k$ of the $\Sc_{2k}$.
These numbers were found by computing the fractional parts
of all $B_{2k}$ for $2k\le10^m$ via a sieve, as in \cite{EW}, sorting them and
counting the unique values.  High enough precision was used
to distinguish close but not equal fractional parts.  

There are a number of problems where the expression $x/(\log x)^\beta$ arises,
and in some of these a secondary factor of $(\log\log x)^c$ arises in the
denominator, sometimes with $c=3/2$ (see \cite{F}) and sometimes with $c=1/2$
(see \cite{BNPT} and \cite{F2}).  We have no compelling reason to suggest such a factor
here, but we've noticed that $F(x)$ has a ratio with 
$f(x):=x/((\log x)^\beta(\log\log x)^{1/2})$ that stays fairly constant. 
 In Table \ref{Ta:countF} we have also recorded
the ratios $R(x)=F(x)/f(x)$.
\begin{table}[ht]
\caption{Number of elements $\le10^m$ of $\Fc$.}
\label{Ta:countF}
\begin{tabular}{|rrrrrr|} \hline
$m$   & 5 & 6 & 7 & 8 & 9 \\ 
$F(10^m)$   & 24662 & 235072 & 2261011 & 21876975 &
212656697 \\ 
$R(10^m)$&.476&.478&.479&.480&.481\\
\hline
\end{tabular}
\end{table}


Table \ref{Ta:countD} shows the growth rate of $\Dc$,
the set of Bernoulli denominators.  These counts were computed with
Mathematica using the criterion that even $d>2$ is in $\Dc$ if and only if
$d=D_{\lambda(d)}$.  The counts were then checked via the sieve procedure
that we used for Table \ref{Ta:countF}.  Note that Theorem \ref{thm:Ddist}
does {\it not} assert that $D(x)/\pi(x)$ tends to a limit or is bounded, but we do know
that these ratios have a positive liminf.

\begin{table}[ht]
\caption{Number of elements $\le10^m$ of $\Dc$.}
\label{Ta:countD}
\begin{tabular}{|rrrrrr|} \hline
$m$ & 5 & 6 & 7 & 8 & 9 \\ 
$D(10^m)$ & 513 & 3649 & 27936 & 226111 & 1893060 \\ 
$D(10^m)/\pi(10^m)$&.053&.046&.042&.039&.037\\
\hline
\end{tabular}
\end{table}

\subsection{A partition of the set of primes}
Given an odd prime $p$, the least $n$ with $p\mid D_n$ is evidently $p-1$.
Let $d_p=D_{p-1}/p$.  For example,
\[
d_3=2,~d_5=6,~d_7=6,~d_{11}=6,~d_{13}=210,~d_{17}=30,~d_{19}=42.
\]
For $d$ even, let
\[
\Pc_d=\{p\hbox{ prime}:d_p=d\},
\]
so that the sets $\Pc_d$ give a partition of the odd primes.
\begin{lemma}
\label{lem:partition}
For each odd prime $p$ we have $\lambda(d_p)\mid p-1$.  Further,
$\lambda(d_p)<p-1$ if and only if $d_p\in\Dc$.
\end{lemma}
\begin{proof}
Note that $d_p$ is the product of those primes $q<p$ where $q-1\mid p-1$.
Thus, $\lambda(d_p)$ is a least common multiple of some of the divisors of 
$p-1$, so we must have $\lambda(d_p)\mid p-1$ and 
$D_{\lambda(d_p)}\mid D_{p-1}$.   Also,
 $d_p\mid D_{\lambda(d_p)}$ since this holds for all squarefree numbers
 larger than 2.  If $\lambda(d_p)=p-1$, then $p\mid D_{\lambda(d_p)}$.
 But $p\nmid d_p$, so $d_p\ne D_{\lambda(d_p)}$ and $d_p\notin\Dc$.
 If $\lambda(d_p)<p-1$, then $D_{\lambda(d_p)}\mid d_p$ since
 $q\mid D_{\lambda(d_p)}$ implies $q-1\mid\lambda(d_p)$, which
 implies $q-1\mid p-1$ and $q<p$, so that $q\mid d_p$.  Hence
 $d_p=D_{\lambda(d_p)}$ and $d_p\in\Dc$.
 \end{proof}
 
 A consequence of Lemma \ref{lem:partition} is that if $d\notin\Dc$,
 then $\#\Pc_d\le 1$.  Indeed, if $p\in\Pc_d$ with $d\notin\Dc$,
 the lemma implies that $\lambda(d)=p-1$, so that $p$ is uniquely
 determined from $d$.  On the other hand, in the next theorem, we
 see that if $d\in\Dc$, $\Pc_d$ is a quite thick set of primes.

\begin{theorem}
\label{thm:dens}
For each $d\in\Dc$ there is a positive constant $c_d$ such that the
relative density of $\Pc_d$ in the set of prime numbers is $c_d$.
\end{theorem}
\begin{proof} (Sketch.)
First note that if $p\in\Pc_d$ then $p\equiv1\pmod{\lambda(d)}$.
 From \cite[Theorem 3]{LPP} it follows that there is an absolute constant $c>0$
 such that for any $3\le z\le x$, the number of primes $p\le x$ such that $p-1$
 has a divisor of  the form $q-1$ with $q$ prime and $z<q<p$ is $O(\pi(x)/(\log z)^c)$.
 This result is completely analogous to \cite[Theorem 2]{EW}.
 We apply this to primes $p\equiv1\pmod{\lambda(d)}$, which comprise a positive
 proportion of all primes by the prime number theorem for residue classes.
 So, it follows from the method in \cite[Theorem 3]{EW} that the set
 \[
 \{p\,\equiv\,1\kern-5pt\pmod{\lambda(d)}:q-1\mid p-1\hbox{ implies }q-1\mid\lambda(d)\hbox{ or }q=p\},
 \]
 where $p,q$ are understood as primes, has a positive relative density $c_d$ in
 the set of all primes.
 Since 
 $\prod_{q-1\mid\lambda(d)}q=d$ by Proposition \ref{prop:lambda},
 for such primes $p$ we have $D_{p-1}=dp$, completing the proof.
 \end{proof}

 We conjecture that $c_6$ is the largest of the densities $c_d$.
 Table \ref{Ta:countdp} has  some counts for $d=6,30,42,66$
 plus fractions of all primes to the same bounds.
 \begin{table}[ht]
\caption{Number of primes $p\le10^k$ with $D_{p-1}=dp$ and fraction of all primes
to $10^k$}
\label{Ta:countdp}
\begin{tabular}{|rrrrrr|} \hline
$d$ & $k=5$ & 6 & 7 & 8 & 9 \\ \hline
6 &1135 & 8772 & 71421& 601804& 5189442\\
&.1183 & .1117 & .1075 & .1045 & .1021 \\
30 & 600& 4312 & 34065& 278709 & 2358192 \\
& .0626 & .0549 & .0513 & .0484 & .0464\\
42 & 480 & 3543 & 27722& 226087 & 1896172\\
& .0500 & .0451 & .0417 & .0392 & .0373 \\
66 & 275& 1933& 14859 & 120565& 1010251 \\ 
& .0287 & .0246 & .0224 & .0209 & .01999\\
 \hline
\end{tabular}
\end{table}
  
  We noticed in tabulating $\Dc$ that there are quite a few more
values of $d\in\Dc$ with $d/6\equiv2\pmod 3$ than with $d/6\equiv1\pmod 3$.
This phenomenon may be due to the robust size of $\Pc_6$ as seen in
Table \ref{Ta:countdp}: every member of $\Pc_6$ when divided by 6 is
$\equiv2\pmod3$.  To be sure, the other cases in Table \ref{Ta:countdp}
count against this trend, but when counting members of $\Dc$ up to
$x$, the $\Pc_6$ members involve primes to $x/6$, while the other
cases involve primes to $x/30$ and smaller.  There may well be other
forces at play here, but this observation may partially demystify the
phenomenon.   We don't know if this imbalance continues
asymptotically since we don't know if $D(x)$ is of order of magnitude
$\pi(x)$.  (Note that the ``$o(1)$" in Theorem \ref{thm:Ddist} may have
some significance.)

 We have seen in the proof of Theorem \ref{thm:dens} that for many primes
 $p$ we have $D_{p-1}=dp$ with $d\in\Dc$.  However, this is not true for
 all odd primes.  For example, note that $D_{12}=210\cdot13$ and $210\notin\Dc$.
 We show that this is uncommon.
 \begin{theorem}
 \label{thm:notinD}
 There is a positive constant $c$ such that the number of primes $p\le x$
 with $D_{p-1}/p\notin\Dc$ is $O(\pi(x)/(\log x)^c)$.  
 \end{theorem}
 \begin{proof}
 Let $p\le x$ be a prime and let $q$ be the largest prime factor of $p-1$.
 Let $B=x^{1/\log\log x}$.  If $q\le B$, that is $p-1$ is a $B$-smooth number,
 we can bound the number of such $p$ by the number of $B$-smooth
 numbers at most $ x$.  By a standard result of de Bruijn, this count is $O_k(x/(\log x)^k)$
for any fixed $k$.  In particular it holds for $k=2$, so we may ignore such primes
and assume that $q>B$.  Next, we again apply \cite[Theorem 3]{LPP} mentioned
in the proof of Theorem \ref{thm:dens}.  We apply this with $z=B$, so the number
of primes $p\le x$ divisible by some shifted prime $r-1$ with $B < r < p$
is negligible.  Thus, we may assume that $p-1$ has no such divisor.  Let $n$ be
the largest $B$-smooth divisor of $p-1$, so that $n\le (p-1)/q<p-1$.  Then
$D_{p-1}=D_np$, so that $D_{p-1}/p\in\Dc$.  This completes the proof.
\end{proof}

\begin{remark}
Similarly, as in Remark \ref{rmk:densum}, we have
\[
\sum_{d\in\Dc}c_d=1.
\]
This follows from Theorem \ref{thm:notinD} and from the fact that for large $B$,
the primes $p$ divisible by some shifted prime $q-1>B$ with $q<p$ are sparse, which follows
from \cite[Theorem 3]{LPP}.
\end{remark}

One might wonder how strong Theorem \ref{thm:notinD} is, or even if there
are infinitely many primes as described in the theorem.  We can prove this
conditionally on the prime $k$-tuples conjecture.  Further,  from the Hardy--Littlewood 
quantitative form of $k$-tuples, we can show there are quite a few of these primes.
\begin{theorem}
\label{thm:hl}
Assuming the prime $k$-tuples conjecture, there are infinitely many primes
$p$ with $D_{p-1}/p\notin\Dc$.  Assuming the quantitative form of this
conjecture due to Hardy and Littlewood, the number of such primes $p\le x$
is $\gg \pi(x)/\log x$.
\end{theorem}
\begin{proof}
Let $q$ be a prime with $q\equiv3\pmod4$, $q>3$, and $p=2q-1$ prime.  Let $d$
be such that $D_{p-1}=dp$.  Suppose that $d=D_n$ for some $n$.  Since
$4\mid p-1$ and $p>5$, we have $5\mid d$, so that $4\mid n$.  
Also $q\mid d$, so $q-1\mid n$.
Hence lcm$(4,q-1)=2(q-1)=p-1\mid n$.  This implies that $p\mid D_n=d$ contradicting
$D_{p-1}=dp$ squarefree.  Thus, $D_{p-1}/p\notin\Dc$.  The prime $k$-tuples conjecture implies there are infinitely
many such $p$, and the quantitative form implies that there are $\gg\pi(x)/\log x$
of them at most $x$.  This completes the proof.
\end{proof}

Table \ref{Ta:countnotinDc} illustrates Theorems \ref{thm:notinD} and
\ref{thm:hl} by showing the number of primes $p$
for which $D_{p-1}/p\notin\Dc$.

\begin{table}[ht]
\caption{Number of $p\le10^m$ with $D_{p-1}/p\notin\Dc$.}
\label{Ta:countnotinDc}
\begin{tabular}{|rrrrrr|} \hline
$m$ & 5 & 6 & 7 & 8 & 9 \\ \hline
Count & 4183 & 34647 & 293117 & 2525121 & 22119959 \\
Count$/\pi(10^m)$ & .436 & .441 & .441 & .438 & .435 \\ \hline
\end{tabular}
\end{table}

A puzzle here:  These counts are holding steady at about 44\% of the primes,
yet Theorem \ref{thm:notinD} says the fraction should decay to 0, though maybe
the decay is slow.  This is similar to the decay of $F(10^m)/10^m$ seen in
Table \ref{Ta:countF}, though it definitely seems more glacial in Table
\ref{Ta:countnotinDc}.

\subsection{Further problems}

The first few Bernoulli denominators all have the form $p-1$ for some
primes $p$: 7, 31, 43, 67, 139, 283, 331.
One might wonder whether there are infinitely many $d\in\Dc$
with $d+1$ composite and infinitely many with $d+1$ prime.
Table \ref{Ta:countDprime} shows the fractions of prime and composite
$d+1$ with $d\in\Dc$ below various powers of 10.  It looks like
the composite cases predominate.  A possible proof: 
It likely follows from the proof of
\cite[Theorem 3]{LP} that
\begin{equation}
\label{eq:lplikely}
\#\{\lambda(p-1):p\le x\}\le \pi(x)/(\log x)^{1+o(1)}
\end{equation}
as $x\to\infty$.  We have 
\begin{align*}
\#\{p\le x:p-1\in\Dc\}&=\#\{p\le x:p-1=D_{\lambda(p-1)}\}\\
&\le\#\{\lambda(p-1):p\le x\},
\end{align*}
so that \eqref{eq:lplikely} would imply that $\#\{p\le x:p-1\in\Dc\}$
is bounded above by $\pi(x)/(\log x)^{1+o(1)}$.  Thus, from Theorem \ref{thm:Ddist},
it would be unusual for $d\in\Dc$ to have $d+1$ prime.

A lower bound of similar quality follows from the strong form of the
Hardy--Littlewood conjecture as in Theorem \ref{thm:hl}.   For this, take primes $r$
with $q=2r+1$ prime and $p=6q+1$ prime.  Then $\lambda(p-1)=q-1=2r$,
and $D_{2r}=6q=p-1$, so $p-1\in\Dc$.  By Hardy--Littlewood, the
number of such $p\le x$ is $\gg\pi(x)/(\log x)^2$.

\begin{table}[ht]
\caption{Number and fraction of composite and prime $d+1\le10^m$ for $d\in\Dc$.}
\label{Ta:countDprime}
\begin{tabular}{|rrrrrrrr|} \hline
$m$ & 3 & 4 & 5 & 6 & 7 & 8 & 9 \\ 
\hline
Composite & 4 & 56 & 361 & 2812 & 22759 & 189894 & 1628333 \\
Fraction & .286 & .667 & .704 & .771 & .815 & .840 & .860 \\
\hline
Prime & 10 & 28 & 152 & 837 & 5177 & 36217 & 264727 \\
Fraction & .714 & .333 & .296 & .229 & .185 & .160 & .140 \\
\hline
\end{tabular}
\end{table}

Let $\psi(n)\to\infty$ arbitrarily slowly.  Then the set
\[
\{n\hbox{ even}:D_n>\psi(n)\}
\]
has asymptotic density 0.  This follows from Remark \ref{rmk:densum}.
On the other hand, there are a fair number of $n$ with $D_n$
large:  we have
\[
\#\{n\hbox{ even}:n\le x,\,D_n>n\}\gg x/\log x.
\]
This follows from Theorem \ref{thm:dens} and the prime number theorem.
In addition, it follows from \cite[Theorem 1]{eps} that
there are positive constants $c,c'$ such that 
$D_n> \exp(n^{c/\log\log n})$ for infinitely many even $n$
and $D_n<\exp(n^{c'/\log\log n})$ always.

In Tables \ref{Ta:count579} and \ref{Ta:count579a} we recorded the second smallest
member of $\Sc_{2k}$ for $2k\in\Fc$ with $2k\le112$.  It would be interesting to study
the distribution of these numbers.

\smallskip
\noindent{\bf Acknowledgment}.
We thank Paul Pollack for reminding us of \cite{LPP}.  And we thank the
referee for a careful reading.



\begin{thebibliography}{99}

\bibitem{BNPT} 
M. Balazard, J.L. Nicolas, C. Pomerance, and G. Tenenbaum, 
Grandes d\'eviations pour certaines fonctions arithm\'etiques,
{\it J. Number Theory} {\bf40} (1992), 146--164.

\bibitem{E48}
P. Erd\H os, On the density of some sequences of integers,
{\it Bull.\ Amer.\ Math.\ Soc.}
{\bf54} (1948), 685--692.

\bibitem{eps}
P. Erd\H os, C. Pomerance, and E. Schmutz,
Carmichael's lambda function, {\it Acta Arith.} {\bf 58} (1991), 363--385.

\bibitem{EW}
P. Erd\H os and S. S. {Wagstaff, Jr.}, The fractional parts of the Bernoulli
numbers, {\it Illinois J.\ Math.}
{\bf24} (1980), 105--112.

\bibitem{F}
K. Ford, The distribution of integers with a divisor in a given interval, 
{\it Ann.\ of Math.} {\bf168} (2008), 367--433.

\bibitem{F2}
K. Ford, Integers divisible by a large shifted prime,
{\it Acta Arith.} {\bf178} (2017), 163--180.

\bibitem{FLP}
K. Ford, F. Luca, and C. Pomerance, The image of Carmichael's $\lambda$-function,
{\it Algebra \& Number Theory} {\bf8} (2014), 2009--2026. 

\bibitem{HR}
H. Halberstam and H.-E. Richert, Sieve Methods, London Mathematical Society Monographs,
No.\ 4. Academic Press [A subsidiary of Harcourt Brace Jovanovich, Publishers], London -- New York, 1974.

\bibitem{HRo}
H. Halberstam and K. F. Roth, Sequences, Volume 1, Oxford University Press,
1966.

\bibitem{LPP}
F. Luca, A. Pizarro-Madariaga, and C. Pomerance, 
On the counting function of irregular primes, 
{\it Indag.\ Math.} {\bf26} (2015), 147--161.

\bibitem{LP}
F. Luca and C. Pomerance,
On the range of Carmichael's universal-exponent function,
{\it Acta Arith.} {\bf 162} (2014), 289--308.

\bibitem{MPP}
N. McNew, P. Pollack and C. Pomerance,
Numbers divisible by a large shifted prime and large torsion subgroups of
CM elliptic curves,
{\it Int. Math. Res. Not.} 2017; doi: 10.1093/imrn/rnw173.

\bibitem{PS}
C. Pomerance and A. S\'ark\"ozy, On homogeneous multiplicative hybrid problems
in number theory, {\it Acta Arith.} {\bf  49} (1988), 291--302.

\bibitem{Sun}
R. Sunseri,
Zeros of $p$-adic $L$-functions and densities related to Bernoulli numbers,
Ph.D. thesis at University of Illinois, Urbana, IL, 1980, available at
{\tt http://homes.cerias.purdue.edu/$\sim$ssw/sunseri.pdf}.

\end{thebibliography}
\end{document}